\documentclass[11pt]{amsart}
\usepackage{amssymb, amsthm, graphicx, graphs, epsfig, verbatim, amscd,
enumerate, stmaryrd, amsmath, psfrag}

\newtheorem{thm}{Theorem}[section]
\newtheorem{cor}[thm]{Corollary}

\newtheorem{prop}[thm]{Proposition}

\theoremstyle{definition}
\newtheorem{defn}[thm]{Definition}

\theoremstyle{remark}

\numberwithin{equation}{section}

\newcommand{\al}{\alpha} 
\newcommand{\be}{\beta} 
\newcommand{\ga}{\gamma} 
\newcommand{\Ga}{\Gamma}

\newcommand{\ep}{\varepsilon}

\newcommand{\ka}{\kappa}

\newcommand{\si}{\sigma} 
\newcommand{\Si}{\Sigma}

\newcommand{\va}{\varphi} 
\newcommand{\csi}{\xi}

\newcommand{\x}{\times}

\newcommand{\GG}{\mathcal G} 
 
\newcommand{\CC}{\mathcal C}
\newcommand{\BB}{\mathcal B}
\renewcommand{\SS}{{\mathfrak {F}}}
\renewcommand{\aa}{{\mathfrak {a}}}
\newcommand{\bb}{{\mathfrak {g}}}

\renewcommand{\sc}{{\mathfrak {sc}}}
\renewcommand{\c}{{\mathfrak {c}}}

\newcommand{\iii}{{\rm III}}

\newcommand{\imm}{{\mathrm {Imm}}}
\newcommand{\DIFF}{{\mathrm {DIFF}}}
\newcommand{\AUT}{{\mathrm {AUT}}}

\newcommand{\ISO}{{\mathrm {ISO}}}

\newcommand{\Z}{\mathbb Z} 
 
\newcommand{\Q}{\mathbb Q} 
\newcommand{\R}{\mathbb R} 
\newcommand{\RP}{{\mathbb R}{P}} 
 
\newcommand{\del}{\partial} 
\newcommand{\co}{\colon\thinspace}  

\newcommand{\plto}[1]
 {\xrightarrow{#1}}

\begin{document}
\mathsurround=1pt 
\title[P-T-Sz type construction for non-positive codimension]{Pontryagin-Thom-Sz\H{u}cs type construction for non-positive codimensional singular maps with prescribed singular fibers}

\address{Faculty of Mathematics, Kyushu University, 6-10-1 Hakozaki, Higashi-ku, Fukuoka 812-8581, Japan}
\email{kalmbold@cs.elte.hu}


\thanks{The author was supported by Canon Foundation in Europe and has been supported by JSPS}




\author{Boldizs\'{a}r Kalm\'{a}r}

\begin{abstract}
We give a Pontryagin-Thom-Sz\H{u}cs type construction for non-positive codimensional singular maps, and
obtain results about cobordism and bordism groups of $-1$ codimensional stable maps with prescribed
singular fibers.
\end{abstract}

\maketitle

\section{Introduction}\label{s:intro}

The purpose of this paper is to show many results
about (co)bordisms of non-positive codimensional singular maps 
(if we have a map $f \co M^m \to P^p$ of an $m$-dimensional manifold 
into a $p$-dimensional manifold, then the {\it codimension} of the map $f$ is the 
integer $p-m$) with prescribed 
singular fibers by
giving a Pontryagin-Thom-Sz\H{u}cs type construction
in an analogous way to the case of positive codimensional singular maps
and by finding analogous theorems and arguments.

The classical Pontryagin-Thom construction is an elementary method to study homotopy groups 
(or stable homotopy groups)
of Thom-spaces of vector bundles via cobordisms of embeddings (resp.\ immersions), or on the 
contrary, study cobordisms via homotopy groups. Sz\H{u}cs extended the classical
Pontryagin-Thom construction for cobordisms of immersions with restricted 
self-intersections \cite{Sz1}, for cobordisms of singular maps with various types of singularities  \cite{Szucs1, Sz2} and in general for cobordisms of positive codimensional singular maps
together with Rim\'anyi \cite{RSz}. Later it turned out that 
Pontryagin-Thom-Sz\H{u}cs type construction can be used effectively to study 
Thom polynomials \cite{Ri} and eliminations of singularities by cobordism \cite{Sz3, Szucs4}. 

In the case of codimension $-1$, we obtained already some results and applications, 
e.g., computing the cobordism groups of fold maps with prescribed singular fibers in low dimensions and eliminating singular fibers
by cobordism \cite{Kal2, Kal10, Kal3, Kal9}.

Recently several techniques have been developed in order to describe cobordisms of singular maps. Ando and Sadykov constructed spectra \cite{An, An6, Sad2} by using the h-principle of Ando \cite{An2, An8}, and Sz\H{u}cs
constructed spectra by using a compression theorem \cite{Szucs4}. The main advantage of the  
Pontryagin-Thom-Sz\H{u}cs type construction is that it also gives the possibility to handle in an elementary way
cobordisms of singular maps with {\it global} restrictions like restrictions about the multiplicities or 
symmetries of singularities or singular or regular fibers.

In Section 2 we give basic definitions, in Section 3 we state our main results, in Section 4 we
give the Pontryagin-Thom-Sz\H{u}cs type construction for non-positive codimensional singular maps,
in Section 5 we reduce the structure groups of $-1$ codimensional singularities to compact subgroups,
in Section 6 we compute some cobordism groups of cusp maps.

The author would like to thank Prof.\ A.\ Sz\H{u}cs for the lectures and discussions.

\section{Preliminaries}\label{s:jelol}

\subsection*{Notations}
In this paper the symbol ``$\amalg$'' denotes the disjoint union.
The symbol $T\csi^k$ denote the 
 the Thom space of the bundle $\csi^k$.
The symbol $\pi_n^s(X)$ (or $\pi_n^s$) denotes the $n$th stable homotopy group of the space $X$ (resp. spheres).
The symbol ``id$_A$'' denotes the identity map of the space $A$.
The symbol $\ep$ denotes a small positive number, and $\Omega_n$
denotes the oriented cobordism group of closed oriented $n$-dimensional manifolds.
All manifolds and maps are smooth of class $C^{\infty}$.

\subsection{Smooth maps}

In this paper, all smooth maps are supposed to be non-positive codimensional, proper and stable. 
We also suppose that a smooth map can be stratified by the singular strata, and the strata are
smooth submanifolds of the target.

\subsection{Bundle structure on a family of a map}

Let $h \co X^l \to \R^k$ be a smooth map of an $l$-dimensional manifold,
where $l \geq k \geq 1$. Let $\GG$ be a subgroup
 of the
{\it automorphism group} $\AUT(h)$ of $h$, i.e., the group of pairs $(\al, \be)$, where  
$\al \co X^l \to X^l$ and $\be \co \R^k \to \R^k$ are diffeomorphisms (suppose that
$\be$ is also a linear transformation), and $h \circ \al = \be \circ h$ holds. 
By the Milnor construction, we can construct the universal ``$h$ bundle''
\[
E\GG \x_{\GG} (h \co X^l \to \R^k) \longrightarrow B\GG,  
\]
i.e., whose ``total space'' is the fiberwise map 
\[
E\GG \x_{\GG} X^l \plto{E\GG \x_{\GG} h} E\GG \x_{\GG} \R^k
\]
denoted by $\chi \co \csi \to \eta$ (the restriction of $\chi$ to any fiber is a map equivalent to $h$).
Now, if we have a map $\va \co W \to B\GG$, then by pull-back, we obtain a 
family 
\begin{center}
\begin{graph}(6,2)
\graphlinecolour{1}\grapharrowtype{2}
\textnode {A}(0.5,1.5){$\va^*\csi$}
\textnode {B}(5.5, 1.5){$\va^*\eta$}
\textnode {C}(3, 0){$W$}
\diredge {A}{B}[\graphlinecolour{0}]
\diredge {B}{C}[\graphlinecolour{0}]
\diredge {A}{C}[\graphlinecolour{0}]
\freetext (3,1.8){$\va^*\chi$}
\end{graph}
\end{center}
of the map $h$ parametrized by the space $W$. Moreover, because of the properties of the Milnor construction,
this family is locally trivial and has structure group $\GG$.

Conversely, if we have a family 
\begin{center}
\begin{graph}(6,2)
\graphlinecolour{1}\grapharrowtype{2}
\textnode {A}(0.5,1.5){$\tilde X$}
\textnode {B}(5.5, 1.5){$\tilde Y$}
\textnode {C}(3, 0){$W$}
\diredge {A}{B}[\graphlinecolour{0}]
\diredge {B}{C}[\graphlinecolour{0}]
\diredge {A}{C}[\graphlinecolour{0}]
\freetext (3,1.8){$\tilde \chi$}
\end{graph}
\end{center}
of the map $h$ parametrized by the space $W$, and this family is locally 
trivial and has structure group $\GG$, then it can be induced from the universal $h$ bundle 
$\chi \longrightarrow B\GG$.

Since in the case of Pontryagin-Thom-Sz\H{u}cs type construction for non-positive codimension
we want to induce families of maps from universal bundles obtained by the Milnor construction similarly to above,
our families in hand should be locally trivial. Moreover, since we want to use transversality in
the total space of the bundle $E\GG \x_{\GG} \R^k \to B\GG$, the topological group $\GG$ 
should be compact or finite dimensional Lie group. In the case of positive codimension, 
see \cite{Jan, RSz, Szucs3, Wa}.
 

\subsection{Fiber-germs, singular fibers and $\tau$-maps}

Let $f \co Q^{q} \to N^{n}$ be a smooth map between smooth manifolds
of dimensions $q$ and $n$ 
respectively, $q \geq n \geq 1$, and $p \in Q^{q}$. 
The {\it fiber-germ} $\SS$ over $f(p)$ is the map germ 
$$f \co (f^{-1}(U_{f(p)}), f^{-1}(f(p))) \to (U_{f(p)}, f(p)),$$ where
$U_{f(p)}$ denotes a small neighbourhood of the point $f(p)$ (see ``singular and regular fibers'' in \cite[Chapter~1]{Sa}).

In this paper, we do not make difference between a fiber-germ $\SS$ and its suspensions
\begin{multline*}
f\x {\mathrm {id}}_{\R^s} \co 
((f\x {\mathrm {id}}_{\R^s})^{-1}(U_{f(p)} \x \R^s), (f\x {\mathrm {id}}_{\R^s})^{-1}(f(p) \x \{0\}) ) \to \\
(U_{f(p)} \x \R^s, f(p) \x \{0\}),
\end{multline*}
$s >0$, although they are not right-left equivalent. E.g., we will always suppose that a 
fiber-germ $\SS$ is not a suspension.

In this paper, instead of a fiber-germ $\SS$ (which is never a suspension by the previous paragraph), 
often we consider a representative
$\si_\SS \co s_\SS \to D_\ep^m$, where
$D_\ep^m$ is a small $m$-dimensional disk around $f(p)$,
$s_\SS$ denotes $f^{-1}(D_\ep^m)$ and $\si_\SS$ denotes
$f|_{f^{-1}(D_\ep^m)}$ (and $\si_\SS$ is not right-left equivalent to a map of the form 
$\tilde \si_\SS \x {\mathrm {id}}_{\R^s}$, $s>0$). Let $\ka(\SS)$ denote $m$.
We will refer to the representative $\si_\SS \co s_\SS \to D_\ep^m$ as the fiber-germ itself.

If $p \in Q^q$ is a singular point of the map $f$, then we call the fiber-germ over $f(p)$ a {\it singular fiber}.
If the fiber-germ has only regular points, then we call the fiber-germ over $f(p)$ a {\it regular fiber}, 
and the representative $\si_\SS \co s_\SS \to D_\ep^m$ is chosen to be
the map $f|_{f^{-1}(f(p))} \co f^{-1}(f(p)) \to f(p)$ (i.e., $\ka(\SS)=m=0$).

\begin{defn}
Let $f \co Q^{q} \to N^{n}$ be a smooth map. Let $\SS$ be a fiber-germ of $f$ and let $S_\SS$ denote
the $\ka(\SS)$ codimensional submanifold of $N^n$ over which the fiber-germs are equivalent to $\SS$. 
We say that the {\it local triviality condition} holds for the family of $\si_\SS \co s_\SS \to D_\ep^{\ka(\SS)}$
parametrized by $S_\SS$ if it is a locally trivial bundle with fiber $\si_\SS$ and structure group 
$\AUT(\si_\SS)$. Furthermore, we say that the {\it compact structure group condition} holds
for the family of $\si_\SS \co s_\SS \to D_\ep^{\ka(\SS)}$
parametrized by $S_\SS$ if its structure group can be reduced to a compact subgroup of $\AUT(\si_\SS)$.
\end{defn}

For example, a family of a regular fiber $\si_\SS \co s_\SS \to D^0$ always satisfies the local triviality condtition, and if
the self-diffeomorphism group of the
manifold $s_\SS$ can be reduced to a compact group (e.g. when $s_\SS$ is one dimensional), then
it satisfies the compact structure group condition as well.

\begin{defn}
We say that a smooth map $f$ is {\it locally trivial} if all of its fiber-germ families satisfy the local 
triviality condition. We say that a locally trivial smooth map $f$ is {\it compact} if all of its fiber-germ families satisfy the compact structure group condition.
\end{defn}

We introduce some notion which are analogues of ``$\tau$-maps'' of the 
positive codimensional case \cite{RSz}. 

\begin{defn}
Let $\si_\SS \co s_\SS \to D_\ep^{\ka(\SS)}$ be a fiber-germ. 
Let $g_1$ be an auto-diffeomorphism of $s_\SS$ and $g_2$ 
be an auto-diffeomorphism of $D_\ep^{\ka(\SS)}$ such that
$g_2 \circ \si_\SS = \si_\SS \circ g_1$.
We call the pair $(g_1, g_2)$ an {\it automorphism} of the fiber-germ $\si_\SS$.
The {\it automorphism group} of $\si_\SS$ consists of this kind of pairs $(g_1, g_2)$.
Let $\AUT(\si_{\SS})$ denote the automorphism group of $\si_\SS$. 
\end{defn}

Hence $\AUT(\si_{\SS})$ is a subgroup of the topological group $\DIFF(s_\SS) \x \DIFF(D_\ep^{\ka(\SS)})$.

\begin{defn}
Let $\si_\SS$ be a fiber-germ and let $\GG$ be a subgroup of the automorphism group of $\si_\SS$. 
We call the pair $(\si_\SS, \GG)$ a {\it global fiber-germ}.
\end{defn}

\begin{defn}
Let $\tau_0$ and $\tau_1$ be two sets of global fiber-germs. $\tau_0 \prec \tau_1$
if for every $(\si_\SS, \GG_\SS) \in \tau_0$ there exists an element 
$(\si_\SS', \GG_\SS') \in \tau_1$ such that $\si_\SS = \si_\SS'$ and 
$\GG_\SS$ is a subgroup of $\GG_\SS'$.
\end{defn}

\begin{defn}
For a locally trivial map $f$ 
let $\tau_f$ be the set of the global fiber-germs $(\si_\SS, \GG_\SS)$
where each $\si_\SS$ is a fiber-germ of $f$ and the corresponding group $\GG_\SS$ is 
the structure group of the bundle of this fiber-germ. 
We call the set $\tau_f$ 
the {\it global fiber-germ set of $f$}.
\end{defn}

\begin{defn}
Let $\tau$ be a set of global fiber-germs. 
We say that the locally trivial map $f$ is a {\it $\tau$-map} if $\tau_f \prec \tau$.
\end{defn}

\subsection{Some $-1$ codimensional singular maps}
 
Let $Q^{n+1}$ and $N^n$ be smooth manifolds of dimensions $n+1$ and $n$ 
respectively. Let $p \in Q^{n+1}$ be a singular point of 
a smooth map $f \co Q^{n+1} \to N^{n}$. The smooth map $f$  has a {\it fold 
singularity} at the singular point $p$, if we can write $f$ in some local coordinates at $p$  
and $f(p)$ in the form 
\[  
f(x_1,\ldots,x_{n+1})=(x_1,\ldots,x_{n-1}, x_n^2 \pm x_{n+1}^2).
\] 
A smooth map $f \co Q^{n+1} \to N^{n}$ is called a {\it fold map}, if $f$ has only 
fold singularities.

The smooth map $f$  has a {\it cusp 
singularity} at the singular point $p$, if we can write $f$ in some local coordinates at $p$  
and $f(p)$ in the form 
\[  
f(x_1,\ldots,x_{n+1})=(x_1,\ldots,x_{n-1}, x_n^3 + x_1x_n - x_{n+1}^2).
\] 
A smooth map $f \co Q^{n+1} \to N^{n}$ is called a {\it cusp map}, if $f$ has only 
fold and cusp singularities.

The possible regular and singular fibers of stable fold and cusp maps for $n \leq 4$ are 
classified in \cite{Lev, Sa, SaYa, Ya}.

We will prove that every $-1$ codimensional stable map is locally trivial and compact.

\subsection{Bordisms and cobordisms of fold maps with prescribed fiber-germs}\label{kob}

In the following, we suppose that the singular maps are locally trivial. 

\begin{defn}[Cobordism]\label{cobdef}
Let $N^n$ be an $n$-dimensional manifold.
Let $\tau$ and $\tau'$  be two sets of global fiber-germs.
Two $\tau$-maps $f_0$ and $f_1$ 
of closed $q$-dimensional manifolds $Q_0^q$ and $Q_1^q$ into $N^n$ are  
{\it $\tau'$-cobordant}, if 
\begin{enumerate}
\item
there exists a $\tau'$-map  
$F \co X^{q+1} \to N^n \times [0,1]$ from a compact $(q+1)$-dimensional 
manifold $X^{q+1}$,
\item
$\del X^{q+1} = Q_0^q \amalg Q_1^q$,
\item
${F |}_{Q_0^q \x [0,\ep)}=f_0 \x
{\mathrm {id}}_{[0,\ep)}$ and ${F |}_{Q_1^q \x (1-\ep,1]}=f_1 \x {\mathrm {id}}_{(1-\ep,1]}$, where $Q_0^q \x [0,\ep)$
 and $Q_1^q \x (1-\ep,1]$ are small collar neighbourhoods of $\del X^{q+1}$ with the
identifications $Q_0^q = Q_0^q \x \{0\}$, $Q_1^q = Q_1^q \x \{1\}$.
\end{enumerate}

We call  the map $F$ a {\it $\tau'$-cobordism} between $f_0$ and $f_1$.
\end{defn} 
This clearly defines an equivalence relation on the set of $\tau$-maps of closed
$q$-dimensional manifolds into $N^n$.

For two sets $\tau$ and $\tau'$ of global fiber-germs  
let us denote the cobordism classes of $\tau$-maps under $\tau'$-cobordisms 
by $$\CC ob_{N,\tau,\tau'}^{}(q,n-q)$$ ($\CC ob_{N,\tau}^{}(q,n-q)$ for $\tau = \tau'$).

When the target manifold $N^n$ has the form $\R^1 \x M^{n_1}$, we can define a commutative semigroup operation on $\CC ob_{N,\tau,\tau'}^{}(q,n-q)$ 
by the far away disjoint union. In the case of $N^n=\R^n$, 
we have an abelian group, which we denote by
${\mathcal Cob}_{\tau,\tau'}^{}(q, n-q)$ (${\mathcal Cob}_{\tau}^{}(q, n-q)$ for $\tau=\tau'$).

\begin{defn}[Bordism]
Let $N_i^n$ $(i=0,1 )$ be two closed oriented $n$-dimensional manifolds.
Let $\tau$ and $\tau'$  be two sets of global fiber-germs.
Two $\tau$-maps $f_i$ $(i=0,1)$ 
of closed $q$-dimensional manifolds $Q_0^q$ and $Q_1^q$ into $N^n_0$ and $N^n_1$, respectively, are  
{\it $\tau'$-bordant}, if 
\begin{enumerate}
\item
there exists a $\tau'$-map  
$F \co X^{q+1} \to Y^{n+1}$ of a compact $(q+1)$-dimensional 
manifold $X^{q+1}$ to a compact oriented $(n+1)$-dimensional manifold $Y^{n+1}$,
\item
$\del X^{q+1} = Q_0^q \amalg Q_1^q$ and
$\del Y^{n+1} = N_0^n \amalg -N_1^n$,
\item
${F |}_{Q_0^q \x [0,\ep)}=f_0 \x
{\mathrm {id}}_{[0,\ep)}$ and ${F |}_{Q_1^q \x (1-\ep,1]}=f_1 \x {\mathrm {id}}_{(1-\ep,1]}$, where $Q_0^q \x [0,\ep)$
 and $Q_1^q \x (1-\ep,1]$ are small collar neighbourhoods of $\del X^{q+1}$ with the
identifications $Q_0^q = Q_0^q \x \{0\}$, $Q_1^q = Q_1^q \x \{1\}$.
\end{enumerate}

We call  the map $F$ a {\it $\tau'$-bordism} between $f_0$ and $f_1$.
\end{defn} 
This clearly defines an equivalence relation on the set of
$\tau$-maps of closed $q$-dimensional manifolds into closed oriented $n$-dimensional manifolds.

For two sets $\tau$ and $\tau'$ of global fiber-germs
let us denote the bordism classes of $\tau$-maps under $\tau'$-bordisms 
by $\BB or_{\tau,\tau'}^{}(q,n-q)$ ($\BB or_{\tau}^{}(q,n-q)$ for $\tau = \tau'$).

We define a commutative group operation on $\BB or_{\tau,\tau'}^{}(q,n-q)$ in the usual way by the disjoint union.

\section{Main theorems}

The following theorem is the non-positive codimensional analogue of \cite{RSz}.

\begin{thm}\label{PTcons}
For every integers $q \geq n \geq 1$ and set of global fiber-germs $\tau$,
there is a Pontryagin-Thom-Sz\H{u}cs type construction for 
locally trivial $\tau$-maps of $q$-dimensional 
manifolds into $n$-dimensional manifolds, i.e.,
\begin{enumerate}[\textup{(}\rm 1\textup{)}]
\item
there exists a universal $\tau$-map 
$\csi_{\tau}^{}  \co U_{\tau} \to \Ga_{\tau}$ such that\footnote{The spaces $U_{\tau}$ and $\Ga_{\tau}$
are not (finite dimensional) manifolds and so $\csi_{\tau}^{}$ is not a smooth map\textup{.}}
for every locally trivial $\tau$-map $g \co Q^{q} \to N^{n}$ there exists a commutative diagram 
\[
\begin{CD}
Q^{q} @>>> U_{\tau} \\
@V g VV @V \csi_{\tau}^{} VV \\
N^{n} @>>> \Ga_{\tau}
\end{CD}
\]
moreover the arising map $N^{n} \to \Ga_{\tau}$ is unique up to homotopy\textup{.} It will be denoted by $\chi_g$\textup{.}
The space $\Ga_{\tau}$ is constructed by gluing together total spaces of vector bundles
corresponding to the possible singular fibers in $\tau$ and their automorphisms\textup{.}
\item
For every $n$-dimensional manifold $N^n$ 
there is a natural map 
\[
\chi_*^N \co \CC ob_{N^n,\tau}^{}(q,n-q) \to [{\dot N}^n, \Ga_{\tau}]
\]
between the set of $\tau$-cobordism classes $\CC ob_{N^n,f}^{}(n+1,-1)$ 
and the set of homotopy classes $[{\dot N}^n, \Ga_{-1}]$\textup{.} 
The map $\chi_*^N$ maps a $\tau$-cobordism class
$[g]$ of a locally trivial map $g$ into the homotopy class of the inducing map $\chi_g \co {\dot N}^{n} \to \Ga_{\tau}$\textup{.}  
\item
If $n,q,\tau$ are such that the $\tau$-maps and their cobordisms are not only locally trivial but also compact, then the 
natural map $\chi_*^N$ is a bijection.
\end{enumerate}
\end{thm}

Similarly, there is a natural map $\BB or_{\tau}^{}(q,n-q) \to \Omega_n(\Ga_{\tau})$, which is a 
bijection if the maps in hand are compact.

\begin{defn}\label{weakcoo}
Let $\si_\SS \co s_\SS \to D_\ep^{\ka(\SS)}$ be a singular fiber.
We say that the singular fiber $\si_\SS$ is {\it weakly coorientable}
if for every automorphism $(g_1, g_2)$ of $\si_\SS$ the diffeomorphism $g_2$
preserves the orientation. (See similar notions for example in \cite{Sa, SaYa}).
\end{defn}

The following theorem can be obtained by using an analogue of the
spectral sequence argument for example in \cite{Szucs2} and Proposition~\ref{isoveges}.
Let $q \geq n \geq 1$, $q-n=1$, and $\tau$ be a set of global fiber-germs with maximal symmetry groups such that 
\begin{enumerate}
\item
every singular fiber has only fold singularities, 
\item
$\tau$ contains all the regular fibers and definite fold singular fibers, and
\item
if $\SS$ is a singular fiber in $\tau$,  then any singular fiber of the form 
$$\SS \cup {\mathrm {some\ regular\ fibers}} \cup {\mathrm {some\ definite\ fold\ singular\ fibers}}$$
is in $\tau$.
\end{enumerate}

\begin{thm}
The rank of the homology group $H_n(\Ga_{\tau}; \Z)$ is equal
to the rank of the $\mathrm n$th homology group of the universal complex (see \cite{Sa, Sa6}) for weakly coorientable singular fibers 
in $\tau$. If $p$ is an odd prime and $p > n$, then
the homology group $H_n(\Ga_{\tau}; \Z)$ does not have any $p$-torsion.\qed
\end{thm}
 
Since $-1$ codimensional fold maps are compact, we obtain the following.
 
\begin{cor}
If $p$ is an odd prime and $p > n$, then
the cobordism group $\CC ob_{\tau}^{}(n+1,-1)$ does not have any $p$-torsion.
\end{cor}

\begin{cor}
The rational bordism group $\BB or_{\tau}^{}(n+1,-1) \otimes \Q$ is isomorphic to 
$\bigoplus_{k+l = n} H_k(\Ga_{\tau};\Q) \otimes \Omega_l$.
\end{cor}


Analogously to \cite{Kal10},
we can apply the results of \cite{SzaSzu, Szucs4} to the Pontryagin-Thom-Sz\H{u}cs type construction for $-1$ codimensional stable maps with prescribed singular fibers as follows.

For a global singular fiber $(\si_\SS \co s_\SS \to D_\ep^{\ka(\SS)}, \GG)$,
let $\imm_N^\GG(n-\ka(\SS), \ka(\SS))$ 
denote the cobordism 
group of immersions of $(n-\ka(\SS))$-dimensional manifolds into $N^{n}$ with
normal bundle induced from the universal bundle
$$E\GG \x_{\GG} D_\ep^{\ka(\SS)}.$$

Let $q \geq n \geq 1$, $q-n=1$, and $\tau$ be a set of global fiber-germs such that 
\begin{enumerate}
\item
every singular fiber has only stable singularities, 
\item
$\tau$ contains all the regular fibers and definite fold singular fibers, and
\item
if $\SS_1$ and $\SS_2$ are fiber-germs in $\tau$, then the fiber-germ $\SS_1 \cup \SS_2$ (with the
appropriate symmetry group) is also in the set $\tau$.
\end{enumerate}
Let $\si_\SS$ be a connected singular fiber in $\tau$, such that
there is no connected singular fiber in $\tau$ whose boundary has the singular fiber $\SS$ (i.e., $\si_\SS$
is a ``top'' singular fiber).
Let $\tau_0$ denote the set obtained from $\tau$ by leaving out the singular fiber $\si_\SS$
and all of its multi singular fibers.

\begin{thm}\label{klasszterfibr}
There is a fibration
\[
\Gamma_{\tau_0} \hookrightarrow
\Gamma_{\tau} \plto{}  
\Gamma(T(E\GG \x_{\GG} D_\ep^{\ka(\SS)})),
\]
where  $T(E\GG \x_{\GG} D_\ep^{\ka(\SS)})$ denotes the Thom-space of the bundle
$E\GG \x_{\GG} D_\ep^{\ka(\SS)}$ and $\Gamma(T(E\GG \x_{\GG} D_\ep^{\ka(\SS)}))$
denotes the classifying space for the immersions of $\ka(\SS)$ codimensional immersions with
normal bundle induced from the  bundle
$E\GG \x_{\GG} D_\ep^{\ka(\SS)}.$\qed
\end{thm}

\begin{cor}
Let $f \co Q^{n+1} \to N^n$ be a $\tau$-map. 
Then $f$ is
$\tau$-cobordant to a $\tau_0$-map
if and only if the immersion  into
$N^n$ corresponding to the $\SS$-family of $f$ is zero in the cobordism group
$\imm^\GG(n-\ka(\SS), \ka(\SS))$.\qed
\end{cor}


\section{The classifying space for $\tau$-maps}\label{klasszter}

\begin{proof}[Proof of Theorem~\ref{PTcons}]

Given a topological group $G$
let $EG \to BG$ denote the Milnor construction for the 
universal $G$-bundle. For a global fiber-germ $(\si_\SS, \GG)$,
let $\csi_{\SS}^\GG$ denote the ``total space'' of the ``bundle'' 
\[
E\GG \x_{\GG} (\si_{\SS} \co s_{\SS} \to D_\ep^{\ka(\SS)}) \to
B\GG,  
\]
i.e., the fiberwise map
\begin{equation} 
E\GG
\x_{\GG} s_{\SS} \to  E\GG \x_{\GG} D_\ep^{\ka(\SS)}. \label{egyenlet}
\end{equation}
Let us define the boundary $\del \csi_{\SS}^\GG$  of  $\csi_{\SS}^\GG$ as the map
\[
E\GG \x_{\GG} (\del{\si_{\SS}} \co \del{s_{\SS}} \to
\del{D_\ep^{\ka(\SS)}}) 
\]
of $\csi_{\SS}^\GG$. This is a smooth map, the
fiber-germs of which have codimension less than $\ka(\SS)$.

If $\si_\SS$ is a regular fiber, then the boundary $\del \csi_{\SS}^\GG$ is the empty set.

We can define a partial ordering on the set of global fiber-germs.
\begin{defn}
We say that the global fiber-germ $\aa = (\si_{\SS_0}, \GG_0)$ is {\it less} than the
global fiber-germ $\bb = (\si_{\SS_1}, \GG_1)$ (denoted by $\aa < \bb$) if
$\ka(\SS_0) < \ka(\SS_1)$,
$\del \csi_{\SS_1}^{\GG_1}$ has a fiber-germ of type $\si_{\SS_0}$
and the structure group of the fiber-germ $\si_{\SS_0}$ in $\del \csi_{\SS_1}^{\GG_1}$
is a subgroup of the structure group $\GG_0$.
\end{defn}

Note that for any group $\GG_1$ the relation $(\si_{\SS_0}, \AUT(\si_{\SS_0})) < (\si_{\SS_1}, \GG_1)$
holds if and only if the boundary $\del \si_{\SS_1}$ of the fiber-germ $\si_{\SS_1}$ 
has a fiber-germ of type $\si_{\SS_0}$.

\begin{defn}
We say that a global fiber-germ set $\tau$ is {\it closed} if for any 
global fiber-germ $(\si_{\SS}, \GG)$ in $\tau$ and for every fiber-germ 
$\si_\SS'$ of the map $\del \csi_{\SS}^{\GG}$ there exists a global fiber-germ $(\si_{\SS}', \GG')$ 
in $\tau$ such that $(\si_{\SS}', \GG') < (\si_{\SS}, \GG)$.
\end{defn}

Note that for a map $f$ the global fiber-germ set $\tau_f$ is closed.

Conversely, for a closed global fiber-germ set $\tau$, we can construct by induction
a map $\csi_\tau$ such that $\tau_{\csi_\tau}$ is equal to $\tau$ as follows.

The following construction of the map $\csi_\tau$ 
is very similar to the constructions of the classifying spaces of Sz\H{u}cs \cite{ Sz1, Sz2}
and the universal map of Rim\'anyi and Sz\H{u}cs \cite{RSz}. 
We suggest to study and  understand 
the construction in \cite{RSz} before reading the following definition.
In our case the transversality and the uniqueness up to homotopy are 
included without mentioning.

\begin{defn}\label{univdef}
Let $\tau$ be a closed global fiber-germ set. We construct the map $\csi_\tau$ by
induction.
Let $\csi_\tau^0$ denote the disjoint union $\amalg_{\si_\SS} E\GG \x_{\GG} \si_\SS$
where $(\si_{\SS}, \GG)$ runs over the global 
regular fibers in $\tau$ (i.e., $\si_\SS$ whith $\ka(\SS) = 0$).

Suppose that we have constructed the map $\csi_\tau^k$, $(k \geq 0)$. Then for every global fiber-germ
$(\si_{\SS}, \GG)$ in $\tau$ with $\ka(\SS) = k+1$, the map $\del \csi_{\SS}^{\GG}$
has only global fiber-germs in $\tau_{\csi_\tau^k}$ because $\tau$ is closed.
Because of the induction hypothesis, we have
a glueing ``map'' $\va_{\SS}^{\GG} \co \del \csi_{\SS}^{\GG} \to \csi_\tau^k$
between the sources and targets of the maps $\del \csi_{\SS}^{\GG}$ and $\csi_\tau^k$,
which -- restricted to each stratum of the map $\del \csi_{\SS}^{\GG}$ -- is a bundle-map
into the ``total spaces'' $\csi_{\SS'}^{\GG'}$ (where $\csi_{\SS'}^{\GG'} \subseteq \csi_\tau^k$) 
for every global fiber-germ $(\si_{\SS'}, \GG')$ with 
$(\si_{\SS'}, \GG') < (\si_{\SS}, \GG)$.
Now, by using these glueing maps, 
we obtain the desired map $\csi_\tau^{k+1}  = \csi_\tau^k \cup_{\va_{\SS}^{\GG}} \csi_{\SS}^{\GG}$, where
$(\si_{\SS}, \GG)$ runs over the set of global fiber-germs in $\tau$ with $\ka(\SS)= k+1$.
Let $\csi_\tau$ denote the limit $\csi_\tau^1 \subseteq \csi_\tau^2 \subseteq \ldots \subseteq \csi_\tau^k \subseteq \ldots$.
\end{defn}

The map $\csi_\tau$ defined in Definition~\ref{univdef} has a target space $\Ga_\tau$, i.e., the target spaces 
\[
\Ga_{\SS}^\GG := E\GG \x_{\GG} D_\ep^{\ka(\SS)} 
\]
of the maps $\csi_{\SS}^\GG$, see (\ref{egyenlet}), glued together.

\begin{prop}
The map  $\csi_\tau \co U_\tau \to \Ga_\tau$  is the Pontryagin-Thom-Sz\H{u}cs type construction
for Theorem~\ref{PTcons}.
\end{prop}

\begin{proof}
The homopoty groups of $\Ga_\tau$ provide us the corresponding 
cobordism groups according to the Pontryagin-Thom-Sz\H{u}cs type
construction as follows.

Let $f \co Q^q \to N^n$ be a  $\tau$-map. 
Then, we obtain a map $\mathfrak S \co \dot N^n \to \Ga_\tau$ as follows.
Let $U_{\SS(f)}$ denote the small tubular neighbourhood of the $f$-image\footnote{In the case of a regular fiber let
$U_{\SS(f)}$ denote only the $f$-image.} of the fiber-germs
of type $\si_\SS$ in $N^n$.
By a map $\mathfrak S_{\SS} \co U_{\SS(f)} \to \Ga_{\SS}^\GG$, we induce 
the $\si_\SS$ bundle  of $f$
from the block $\Ga_{\SS}^\GG$ of $\Ga_\tau$
for every global fiber-germ $(\si_{\SS}, \GG')$ in $\tau_f$, such that
the regular neighbourhood $U_{\SS(f)}$ is induced from the target 
 $E\GG \x_{\GG} D_\ep^{\ka(\SS)}$
of the map $\csi_{\SS}^{\GG}$,
and the preimage $f^{-1}(U_{\SS(f)})$ is induced from the source $E\GG \x_{\GG} s_{\SS}$  of the map $\csi_{\SS}^{\GG}$. Then,
we glue together
these maps $\mathfrak S_{\SS} \co U_{\SS(f)} \to \Ga_\tau$ 
into a map $\mathfrak S \co N^n \to \Ga_\tau$ by the
generalized Pontryagin-Thom construction \cite{RSz}.

Conversely, let $\mathfrak S \co \dot N^n \to \Ga_\tau$ be an element of
$[\dot N^n, \Ga_\tau]$.  
If the structure groups of the global fiber-germs in $\tau$ are compact, 
then we can construct a $\tau$-map of a closed ${q}$-dimensional manifold into $N^n$
 by glueing together the  $\si_\SS$ bundles 
$\mathfrak S^*(\csi_{\SS}^{\GG}) \to \mathfrak S^{-1}({\mathrm {zero\ section\ of\ }}\Ga_{\SS}^\GG)$
induced from the $\si_\SS$ bundles $\csi_\SS^\GG \to B\GG$.

This can be applied for homotopies and cobordisms as well. Details are left to the reader.
\end{proof}

This completes the proof of Theorem~\ref{PTcons}.
\end{proof}

\section{Bundle structure on $-1$ codimensional stable maps}

Let $f \co Q^{n+1} \to N^n$ be a stable map in general position. 
The map $f$ can be considered as ``locally trivial bundles of singularities'' glued together. 
More precisely, we have the following theorems.
The first is an analogue of \cite{Szucs3} while the second is an analogue of 
the positive codimensional structure group reduction \cite{Jan, RSz, Wa}.

For each singular fiber $\si_\SS$ of the map $f$ 
let $S_{\SS}$ denote the submanifold in $N^n$
which is the $f$-image of singular fibers of type $\si_\SS$. Note that
$S_{\SS}$ is an $(n - \ka(\SS))$-dimensional submanifold.

Let $P_{\SS}$ denote the total space of the disk bundle 
associated with the normal bundle of the submanifold $S_{\SS}$. The manifold $P_{\SS}$ is embedded into $N^n$ in a natural way, onto a regular neighbourhood of $S_{\SS}$. 
Hence, we have the projection $\pi_{\SS} \co P_{\SS} \to S_{\SS}$,
and we also have the projection $f^{-1}(P_{\SS}) \to S_{\SS}$ defined by the map $\pi_{\SS} \circ f'$,
where $f'=f|_{f^{-1}(P_{\SS})}$.
Therefore, we have the commutative diagram
\begin{center}
\begin{graph}(6,2)
\graphlinecolour{1}\grapharrowtype{2}
\textnode {A}(1,1.5){$f^{-1}(P_{\SS})$}
\textnode {B}(3, 0){$S_{\SS}$}
\textnode {C}(5, 1.5){$P_{\SS}$}
\diredge {A}{B}[\graphlinecolour{0}]
\diredge {C}{B}[\graphlinecolour{0}]
\diredge {A}{C}[\graphlinecolour{0}]
\freetext (1.2, 0.6){$\pi_{\SS} \circ f'$}
\freetext (3,1.8){$f'$}
\freetext (4.6,0.6){$\pi_{\SS}$}
\end{graph}
\end{center}
which gives us a $\si_\SS$-family denoted by $\csi_\SS(f)$, i.e., the ``total space'' of $\csi_\SS(f)$ is the fiberwise map
$f' \co f^{-1}(P_{\SS}) \to P_{\SS}$ between the total spaces of the bundles 
$\pi_{\SS} \circ f' \co f^{-1}(P_{\SS}) \to S_{\SS}$ with fiber $s_\SS$
and $\pi_{\SS} \co P_{\SS} \to S_{\SS}$ with fiber $J^{\ka(\SS)}$, 
the ``base space'' of $\csi_\SS(f)$ is $S_{\SS}$,
and the ``fiber'' of $\csi_\SS(f)$ 
is right-left equivalent to the fiber-singularity $\si_\SS$.

Then the following theorem can be proved by an argument similar to that in \cite{Szucs3}.

\begin{thm}\label{szingnyalloktriv}
Let $\si_\SS$ be a singular fiber of the map $f$.
Then the family $\csi_\SS(f)$ is a locally trivial bundle over 
$S_{\SS}$ with the singular fiber $\si_\SS$ as fiber,
and with structure group $\AUT(\si_{\SS})$.\qed
\end{thm}
\begin{proof}
Since the singularities are stable, there is an open covering of $S_\SS$ such that over each open set the family is trivial. 
By doing the same process as in \cite{Szucs3}, we recieve the result.
\end{proof}

\begin{thm}\label{strcsopveges}
Let $\si_\SS$ be a singular fiber of the map $f$.
Then the structure group of the bundle $\csi_\SS(f)$ can be reduced to a compact group.
\end{thm}

\begin{proof}
Let $N(S)$ denote the neighbourhood of the singular set in $f^{-1}(P_{\SS})$. The
restriction 
\begin{center}
\begin{graph}(6,2)
\graphlinecolour{1}\grapharrowtype{2}
\textnode {A}(1,1.5){$N(S)$}
\textnode {B}(3, 0){$S_{\SS}$}
\textnode {C}(5, 1.5){$P_{\SS}$}
\diredge {A}{B}[\graphlinecolour{0}]
\diredge {C}{B}[\graphlinecolour{0}]
\diredge {A}{C}[\graphlinecolour{0}]
\freetext (1.2, 0.6){$\pi_{\SS} \circ f'$}
\freetext (3,1.8){$f'$}
\freetext (4.6,0.6){$\pi_{\SS}$}
\end{graph}
\end{center}
of the bundle  $\csi_\SS(f)$ is a bundle $\nu_\SS(f)$ of a multi-germ around finite number of isolated points, hence
its structure group can be reduced to a maximal compact subgroup by \cite{Jan, Wa}.

The restriction of the map $f$ to the complementer $f^{-1}(P_{\SS}) - N(S)$ is a submersion with one
dimensional manifolds as fibers. Hence, clearly we can put a Riemannian metric on 
$f^{-1}(P_{\SS}) - N(S)$, which is invariant under the compact structure group of the bundle 
$\nu_\SS(f)$ on $(f^{-1}(P_{\SS}) - N(S)) \cap cl(N(S))$.
Therefore the structure group 
$\AUT(\si_{\SS})$
can be reduced to a compact group. 
\end{proof}

Note that we showed a symmetry group, which keeps fixed a Riemannian metric on $s_\SS$. Let
us denote this group by $\ISO(\si_\SS)$. In \cite{Kal10}, we obtained the following. 

\begin{prop}\label{isoveges}
For a singular fiber $\si_\SS$ with only indefinite singular points and no circle components,
the group $\ISO(\si_\SS)$ is finite, and if $p$ is an odd prime with $p > \ka(\SS)$, then its order 
$|\ISO(\si_\SS)|$ cannot be divided by $p$.
\end{prop}

\section{Cusp cobordisms in low dimensions}

By Theorems~\ref{PTcons} and \ref{klasszterfibr} and the previous section, we can compute some 
cobordism groups of $-1$ codimensional cusp maps. Let us denote the set of all global regular fibers and singular fibers with only fold or cusp singularities and 
with maximal structure groups by ``$\c$'', hence  
$\CC ob_\c^{}(n+1,-1)$ denotes the cobordism group of cusp maps of closed $(n+1)$-dimensional 
manifolds into $\R^n$. By introducing orientations on the source manifolds of maps and cobordisms in Definition~\ref{cobdef}, 
we obtain the notion of oriented cobordism group of cusp maps, denoted by $\CC ob_\c^{O}(n+1,-1)$.
Our theorems clearly work for oriented cobordisms as well.

\begin{prop}
The cobordism group $\CC ob_\c^{O}(3,-1)$ is isomorphic to $\Z_2$. 
The homomorphism $\ga \co \CC ob_\c^{O}(3,-1) \to \imm(1,1)$, which
maps a cobordism class $[f]$ into the cobordism class of the immersion of the singular 
set\footnote{Cusp maps restricted to their singular sets are not immersions, but one can consider them
as immersions by ``smoothing'' the cusps.} of $f$, is an isomorphism.
\end{prop}
\begin{proof}
By Theorems~\ref{PTcons} and \ref{klasszterfibr}, we have a homotopy exact sequence for cobordism
groups, and the homomorphism $\ga$ clearly does not vanish, since it can be non-zero 
on fold maps \cite{Kal2}. 
\end{proof}

\begin{prop}
The rank of the cobordism group $\CC ob_\c^{O}(4,-1)$ is $2$. The two $\Z$ summands are generated
by the algebraic number of the singular fiber of type $\iii^8$ and the algebraic number
of the singular fiber of type $\iii^b$ (for the notations of the singular fibers, see \cite{Sa}).
\end{prop}
\begin{proof}
By \cite{Kal9} and the exact sequence from Theorem~\ref{klasszterfibr}, we obtain the result.
\end{proof}

Let us denote the set of all the global regular fibers and singular fibers with only cusp or fold singularities which have
at most one singular point in each of their connected components and 
with maximal structure groups by ``$\sc$'', hence  
$\CC ob_\sc^{O}(n+1,-1)$ denotes the oriented cobordism group of simple cusp maps of closed $(n+1)$-dimensional 
manifolds into $\R^n$. 

\begin{prop}
The simple cusp cobordism group $\CC ob_\sc^{O}(n+1,-1)$ is isomorphic to $\pi_{n-1}^s(\RP^{\infty})$.
\end{prop}
\begin{proof}
By \cite{Kal3} and easy geometric constructions the summand $\pi_{n-1}^s$ of the oriented simple fold cobordism group goes to zero under the natural homomorphism into the group $\CC ob_\sc^{O}(n+1,-1)$.
But by the exact sequence from Theorem~\ref{klasszterfibr}, the kernel of this natural homomorphism cannot be bigger.
\end{proof}

\end{document}